\documentclass[a4paper,11pt]{article}
\usepackage[latin1]{inputenc}
\usepackage[francais,english]{babel}    % Pour du fran\c{c}ais
\usepackage{amssymb}
\usepackage{amsmath}
\usepackage{amsthm}
\usepackage{multicol}
\usepackage{graphicx}
\usepackage{color}
\usepackage[T1]{fontenc}
\usepackage{yfonts}
\usepackage{placeins}
\usepackage{array}
\usepackage{dsfont}
\usepackage{stmaryrd}
\usepackage{verbatim}
\usepackage{caption}
\usepackage{transparent}
\usepackage{bbold}
\usepackage{hyperref}
\usepackage{enumitem}
\usepackage{tikz}
\usepackage{tikz-3dplot}

\newtheorem{thm}{Theorem}[section]
\newcommand{\bthm}{\begin{thm}}
\newcommand{\ethm}{\end{thm}}

\newtheorem{thmi}{Theorem}
\newcommand{\bthmi}{\begin{thmi}}
\newcommand{\ethmi}{\end{thmi}}

\newtheorem{cori}[thmi]{Corollary}
\newcommand{\bcori}{\begin{cori}}
\newcommand{\ecori}{\end{cori}}

\newtheorem{mthm}{Theorem}
\newcommand{\bmthm}{\begin{mthm}}
\newcommand{\emthm}{\end{mthm}}

\newtheorem{mcor}[mthm]{Corollary}
\newcommand{\bmcor}{\begin{mcor}}
\newcommand{\emcor}{\end{mcor}}
\newtheorem{mconj}[mthm]{Conjecture}
\newcommand{\bmconj}{\begin{mconj}}
\newcommand{\emconj}{\end{mconj}}

\newtheorem*{conj}{Conjecture}
\newcommand{\bconj}{\begin{conj}}
\newcommand{\econj}{\end{conj}}

\newtheorem*{question}{Question}
\newcommand{\bq}{\begin{question}}
\newcommand{\eq}{\end{question}}

\newtheorem*{thn}{Theorem}
\newcommand{\bthn}{\begin{thn}}
\newcommand{\ethn}{\end{thn}}

\newtheorem{exo}{Exercise}
\newcommand{\bex}{\begin{exo}}
\newcommand{\eex}{\end{exo}}

\newtheorem{sol}{Solution}
\newcommand{\bsol}{\begin{sol}}
\newcommand{\esol}{\end{sol}}

\newtheorem{pro}[thm]{Proposition}
\newcommand{\bpro}{\begin{pro}}
\newcommand{\epro}{\end{pro}}

\newtheorem{cor}[thm]{Corollary}
\newcommand{\bcor}{\begin{cor}}
\newcommand{\ecor}{\end{cor}}

\newtheorem{lem}[thm]{Lemma}
\newcommand{\blem}{\begin{lem}}
\newcommand{\elem}{\end{lem}}

\theoremstyle{definition}

\newtheorem{defi}[thm]{Definition}
\newcommand{\bdf}{\begin{defi}}
\newcommand{\edf}{\end{defi}}

\newtheorem*{defis}{Definition}
\newcommand{\bdfs}{\begin{defis}}
\newcommand{\edfs}{\end{defis}}

\newtheorem*{rmk}{Remark}
\newcommand{\brk}{\begin{rmk} \upshape}
\newcommand{\erk}{\end{rmk}}

\newtheorem*{rmks}{Remarks}
\newcommand{\brks}{\begin{rmks} \upshape}
\newcommand{\erks}{\end{rmks}}

\newtheorem*{exe}{Example}
\newcommand{\bexe}{\begin{exe} \upshape}
\newcommand{\eexe}{\end{exe}}

\newtheorem*{exes}{Examples}
\newcommand{\bexes}{\begin{exes} \upshape}
\newcommand{\eexes}{\end{exes}}

\newtheorem*{pre}{Proof}
\newcommand{\bp}{\begin{pre} \upshape}
\newcommand{\ep}{\hfill \qed \end{pre}}
\newcommand{\epp}{\end{pre}}

\newcommand{\beq}{\begin{eqnarray*}}
\newcommand{\eeq}{\end{eqnarray*}}

\newcommand{\beqn}{\begin{equation}}
\newcommand{\eeqn}{\end{equation}}

\newcommand{\ben}{\begin{enumerate}}
\newcommand{\een}{\end{enumerate}}
\newcommand{\bit}{\begin{itemize} \renewcommand{\labelitemi}{$\bullet$} \renewcommand{\labelitemii}{$\star$}}
\newcommand{\eit}{\end{itemize}}

\newcommand{\bfg}{
\begin{figure}[H]
\begin{center}}
\newcommand{\efg}{
\end{center}
\end{figure}
\FloatBarrier}

\newcolumntype{M}[1]{>{\raggedright}m{#1}}

\newcommand{\R}{\mathbb{R}}

\newcommand{\Z}{\mathbb{Z}}

\newcommand{\bs}{\symbol{92}}

\newcommand{\lk}{\operatorname{lk}}

\newcommand{\eps}{\varepsilon}
\newcommand{\st}{\, | \,}

\newcommand{\f}{\frac}

\renewcommand{\geq}{\geqslant}
\renewcommand{\leq}{\leqslant}

\renewcommand{\>}{\rangle}

\newcommand{\mk}{\medskip}

\newcommand{\sign}{\begin{flushright}
Thomas Haettel \\
IMAG, Univ Montpellier, CNRS, France \\
thomas.haettel@umontpellier.fr
\end{flushright}}

\makeatletter
\def\Ddots{\mathinner{\mkern1mu\raise\p@
\vbox{\kern7\p@\hbox{.}}\mkern2mu
\raise4\p@\hbox{.}\mkern2mu\raise7\p@\hbox{.}\mkern1mu}}
\makeatother

\makeatletter
\def\maketitles{%
  \null
  \thispagestyle{empty}%
  \vfill
  \begin{center}\leavevmode
    \normalfont
    {\LARGE \@title\par}%
    \vskip 1.2cm
    {\large \@author\par}%
    \vskip 1.2cm
    {\large \@subtitle\par}%
    \vskip 0.8cm
    {\large \@date\par}%
  \end{center}%
  \vfill
  \null
  \cleardoublepage
  }

\def\date#1{\def\@date{#1}}
\def\author#1{\def\@author{#1}}
\def\title#1{\def\@title{#1}}
\def\subtitle#1{\def\@subtitle{#1}}
\makeatother

\usetikzlibrary{arrows}

\title{XXL type Artin groups are CAT(0) and acylindrically hyperbolic}
\author{Thomas Haettel}
\date{\today}

\begin{document}

\selectlanguage{english}

\maketitle

\begin{center}
\begin{minipage}{0.8\textwidth}
\textsc{Abstract.} We describe a simple locally CAT(0) classifying space for XXL type Artin groups (with all labels at least $5$). Furthermore, when the Artin group is not dihedral, we describe a rank $1$ periodic geodesic, thus proving that XXL type Artin groups are acylindrically hyperbolic. Together with Property RD proved by Ciobanu, Holt and Rees, the CAT(0) property implies the Baum-Connes conjecture for all XXL type Artin groups.

\end{minipage}
\end{center}

\let\thefootnote\relax

\section*{Introduction}

Artin-Tits groups are natural combinatorial generalizations of Artin's braid groups. For every finite simple graph $\Gamma$ with vertex set $S$ and with edges labeled by some integer in $\{2,3,\dots\}$, one associates the Artin-Tits group $A(\Gamma)$ with the following presentation:
$$A(\Gamma) = \<S \st \forall \{s,t\} \in \Gamma^{(1)}, w_m(s,t)=w_m(t,s) \mbox{ if the edge $\{s,t\}$ is labeled $m$}.\>,$$
where $w_m(s,t)$ is the word $stst \dots$ of length $m$. Note that when $m=2$, then $s$ and $t$ commute, and when $m=3$, then $s$ and $t$ satisfy the classical braid relation $sts=tst$.

\mk

Also note that when adding the relation $s^2=1$ for every $s \in S$, one obtains the Coxeter group $W(\Gamma)$ associated to $\Gamma$. Most results about Artin-Tits groups only concern particular classes, which we recall now. The Artin group $A(\Gamma)$ is called:
\bit
\item of \emph{large type} if all labels are greater or equal to $3$,
\item of \emph{extra large type} if all labels are greater or equal to $4$,
\item of \emph{extra extra large type} (XXL) if all labels are greater or equal to $5$,
\item \emph{right-angled} if all labels are equal to $2$,
\item \emph{spherical} if $W(\Gamma)$ is finite, and
\item of \emph{type FC} if every complete subgraph of $\Gamma$ spans a spherical Artin subgroup.
\eit

The \emph{rank} of an Artin-Tits group $A(\Gamma)$ is the number of vertices of $\Gamma$. The \emph{dimension} of an Artin-Tits group $A(\Gamma)$ is the largest rank of a spherical Artin subgroup. In particular, every large type Artin group has dimension at most $2$.

\mk

Many geometric questions are still open for general Artin groups (see~\cite{charney_problems} and \cite{mccammond_mysterious}). In particular, Charney asks the following question, to which we believe the answer is positive:

\bmconj \label{mainconj:CAT0}
Every Artin-Tits group is CAT(0), i.e. acts properly and cocompactly on a CAT(0) metric space.
\emconj

This conjecture has been proved for the following classes of Artin groups:
\ben
\item Right-angled Artin groups (see~\cite{charney_davis_salvetti}).
\item Some classes of $2$-dimensional Artin groups (see~\cite{brady_noel_crisp}, \cite{brady_mccammond_artinthree}).
\item Artin groups of finite type with three generators (see~\cite{brady_artin_three}).
\item $3$-dimensional Artin groups of type FC (see~\cite{bell_artin}).
\item The $n$-strand braid group for $n \leq 6$ (see~\cite{brady_mccammond}, \cite{b6}).
\item The spherical Artin group of type $B_4$ (see~\cite{brady_mccammond}).
\een

Since the classes of $2$-dimensional Artin groups studied by Brady and McCammond in (\cite{brady_mccammond_artinthree}) and the extra extra large type Artin groups we are studying in this article have a large intersection, we will state their results more precisely.

\begin{figure}
\begin{center}
\begin{tikzpicture}
\def \p {0.05}
\draw[fill] (-30:1) circle (\p) node(a) {};
\draw[fill] (210:1) circle (\p) node(b) {};
\draw[fill] (90:1) circle (\p) node(c) {};

\draw [->,>=triangle 60] (b) edge (c) (b) edge (a) (c) edge (a);

\draw[fill] (3,-0.5) circle (\p) node(A) {};
\draw[fill] (3+1.5,-0.5) circle (\p) node(B) {};
\draw[fill] (3+1.5,-0.5+1.5) circle (\p) node(C) {};
\draw[fill] (3,-0.5+1.5) circle (\p) node(D) {};

\draw [->,>=triangle 60] (B) edge (A) (B) edge (C) (D) edge (A) (D) edge (C);

\end{tikzpicture}
\end{center}
\caption{The two forbidden subgraphs in Brady and McCammond's result.}
\label{fig:BM}
\end{figure}
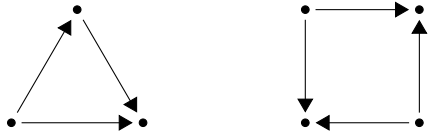

\bthn[Brady and McCammond \cite{brady_mccammond_artinthree}]
Let $A(\Gamma)$ be an Artin group such that one of the following holds:
\bit
\item $|S|=3$ and all labels are greater or equal to $3$.
\item $\Gamma$ contains no triangles.
\item All labels are greater or equal to $3$, and there is a way of orienting the edges of $\Gamma$ so that neither of the graphs in Figure~\ref{fig:BM} appear as subgraphs.
\eit
Then $A(\Gamma)$ is CAT(0).
\ethn

\mk

On the other hand, concerning the cubical world, very few Artin groups have proper and cocompact actions on CAT(0) cube complexes (see~\cite{haettel_artin_cubic}). Nevertheless, the question whether all Artin groups act properly on a CAT(0) cube complex, or more generally have the Haagerup property, is still open.

\mk

Concerning variations on the notion of nonpositive curvature, Bestvina defined a geometric action of Artin groups of spherical Artin on a simplicial complex with some nonpositive curvature features (see~\cite{bestvina_artin}). More recently, Huang and Osajda proved (see~\cite{huang_osajda}) that every Artin group of almost large type (a class including all Artin groups of large type) act properly and cocompactly on systolic complexes, which are a combinatorial variation of nonpositive curvature. They also proved (see~\cite{huang_osajda_helly}) that every Artin group of type FC acts geometrically on a Helly graph, which give rise to classifying spaces with convex geodesic bicombings.

\mk

Another variation on the notion of nonnegative curvature is the acylindrical hyperbolicity (see~\cite{osin_review} for a survey). A group $G$ is called acylindrically hyperbolic if it admits an acylindrical action on some hyperbolic space $X$ (and is not virtually cyclic), i.e. for every $\eps>0$, there exist $N,R \geq 0$ such that, for every $x,y \in X$ at distance at least $R$, we have
$$ |\{g \in G, d(x,g \cdot x) \leq \eps \mbox{ and } d(y,g \cdot y) \leq \eps\}| \leq N.$$
In most cases, it is much easier to find an action on some hyperbolic space with one element satisfying the WPD condition (see~\cite{bestvina_fujiwara}), and then according to Osin (see~\cite{osin}) there exists an acylindrical action on some other hyperbolic space. Concerning Artin-Tits groups, Charney and Morris-Wright (see~\cite{charney_morris_wright}) ask the following question, to which we believe the answer is positive:

\bmconj \label{mainconj:ah}
For every Artin-Tits group $A$, the central quotient $A/Z(A)$ is acylindrically hyperbolic.
\emconj

This conjecture has been proved for the following classes of Artin groups:
\ben
\item Right-angled Artin groups (see~\cite{caprace_sageev}).
\item Braid groups, seen as mapping class groups (see~\cite{masur_minsky} and \cite{bowditch_tight}).
\item Artin-Tits groups of spherical type (see~\cite{calvez_wiest_2}).
\item Artin-Tits groups of type FC such that $\Gamma$ has diameter at least $3$ (see~\cite{chatterji_martin}).
\item Artin-Tits groups such that $\Gamma$ does not decompose as a join of two subgraphs (see~\cite{charney_morris_wright}).
\een

The purpose of this article is to define a new geometric model for XXL type Artin groups.

\bmthm \label{thm:main}
Every XXL type Artin group is the fundamental group of a compact locally CAT(0) $3$-dimensional piecewise Euclidean complex. Furthermore, if the rank of the Artin group is at least $3$, then some element acts as a rank $1$ isometry.
\emthm

An isometry of a CAT(0) space is called \emph{rank 1} if some axis does not bound a flat half-plane. An interesting consequence, due to Sisto (see~\cite{sisto_contracting}), is that if a group $G$ acts properly on a proper CAT(0) space such that some element has rank $1$, then $G$ is either virtually cyclic or acylindrically hyperbolic.

Also note that if $A$ has rank $2$, then $A$ is virtually a direct product of $\Z$ and of a free group, so $A$ is not acylindrically hyperbolic, but its geometry is well understood. In particular, the central quotient $A/Z(A)$ is virtually free and thus acylindrically hyperbolic.

We can deduce the following consequence, regarding the two main conjectures.

\bmcor
Every XXL type Artin group of rank at least $3$ is acylindrically hyperbolic. In particular, Conjecture~\ref{mainconj:CAT0} and Conjecture~\ref{mainconj:ah} hold for all XXL type Artin groups.
\emcor

Note that the class of XXL type Artin groups is not contained in the classes studied by Brady and McCammond in (\cite{brady_mccammond_artinthree}), by Martin and Chatterji (\cite{chatterji_martin}) or by Charney and Morris-Wright (\cite{charney_morris_wright}). For instance, if $\Gamma$ is a complete graph on at least $4$ vertices, with labels at least $5$, then none of the previous results apply.

\mk

Many consequences of being CAT(0) are already consequences of being systolic, and as such are consequences of Huang and Osajda's result (see~\cite{huang_osajda}). For instance, the Novikov conjecture, the fact that centralizers virtually split, the quadratic Dehn function. Let us list a few general consequences of being CAT(0) and acylindrically hyperbolic, which are new for XXL type Artin groups.

\bmcor
Let $A$ be an XXL type Artin group.
\bit
\item $A$ satisfies the $K$-theoretic and $L$-theoretic Farrell-Jones conjectures (see~\cite{bartels_luck} and \cite{wegner_farrell_jones}).
\item $A$ is SQ-universal, i.e. every countable group embeds in a quotient of $A$ (see~\cite{osin}).
\item If $V=\R$ or $V=\ell^p(A)$, for $p \in [1,\infty)$, then $H^2_b(A,V)$ is infinite-dimensional (see~\cite{osin}).
\item $A$ has a free normal subgroup (see~\cite{osin}).
\item $A$ has Property $P_{naive}$: for any finite subset $F \subset A \bs \{1\}$, there exists $g \in A$ such that for all $f \in F$, the group $\<f,g\>$ is freely generated by $\{f,g\}$ (see~\cite{abbott_dahmani}).
\item $A$ is not inner amenable (see~\cite{dahmani_guirardel_osin}).
\item The reduced $C^*$-algebra of $A$ is simple (see~\cite{dahmani_guirardel_osin}).
\eit
\emcor

Ciobanu, Holt and Rees proved (see~\cite{ciobanu_holt_rees}) that every extra large type Artin group satisfies the Rapid Decay Property. According to Lafforgue (see~\cite{lafforgue_baum_connes}), the property RD together with the CAT(0) property imply the Baum-Connes conjecture, so we can state the following.

\bmcor
Every XXL type Artin group satisfies the Baum-Connes conjecture.
\emcor

\textbf{Acknowledgments:} The author would like to thank warmly Chris Cashen for discussions and an invitation to the University of Vienna, where part of this work was initiated. The author would like to thank Anthony Genevois and Damian Osajda for many insightful comments. The author would also like to thank the anonymous referee for comments improving the exposition.
 
\section{The case of dihedral Artin groups}

We start by decribing a very simple nonpositively curved metric model for dihedral Artin groups, which we will use in the sequel as building blocks. If $m \geq 2$, let us denote the dihedral Artin group by $I_2(m)=\<a,b \st w_m(a,b)=w_m(b,a)\>$.

\blem
For every $m \geq 5$, there exists a compact, locally CAT(0), $3$-dimensional piecewise Euclidean complex $X_m$ and $x_0 \in X_m$ with $\pi_1(X_m,x_0) \simeq I_2(m)=\<a,b \st w_m(a,b)=w_m(b,a)\>$. There exist locally geodesic oriented loops $X_m^a$, $X_m^b$ of length $1$ through $x_0$ such that $\pi_1(X_m^a,x_0) = \<a\>$ and $\pi_1(X_m^b,x_0) = \<b\>$. Let $a^+,a^- \in \lk_{x_0}(X_m)$ denote the images in the link of $x_0$ of the positive and negative sides of the loop $X_m^a$, and similarly $b^+,b^- \in \lk_{x_0}(X_m)$ for $X_m^b$. We have furthermore:
\bit
\item $X_m^a \cap X_m^b = \{x_0\}$,
\item $\sphericalangle_{x_0}(a^+,b^+)=\sphericalangle_{x_0}(a^-,b^-) > \frac{4\pi}{5}$.
\item $\sphericalangle_{x_0}(a^+,b^-)=\sphericalangle_{x_0}(a^-,b^+) > \frac{3\pi}{5}$.
\eit
In addition, if $m \geq 6$, we have $\sphericalangle_{x_0}(a^+,b^-)=\sphericalangle_{x_0}(a^-,b^+) > \frac{2\pi}{3}$.
\elem

\bp
Fix $\alpha \in (0,\tan (\f{\pi}{10}))$.

\mk

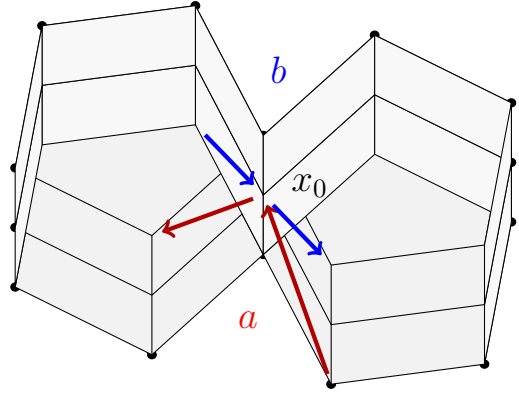
\begin{figure}
\begin{center}
\tdplotsetmaincoords{20}{10}
\begin{tikzpicture}[tdplot_main_coords,scale=1.1547]
\def \p {0.05}
\def \r {1.5}
\def \h {2}
\def \op {1}
\def \oop {1}
\def \gris {black!5}
\def \griss {black!3}

\draw[fill] (\r,0,0) circle (\p) node(11) {};
\draw[fill] (0.309*\r,0.951*\r,0) circle (\p) node(12) {};
\draw[fill] (-0.809*\r,0.588*\r,0) circle (\p) node(13) {};
\draw[fill] (-0.809*\r,-0.588*\r,0) circle (\p) node(14) {};
\draw[fill] (0.309*\r,-0.951*\r,0) circle (\p) node(15) {};

\draw[fill] (\r,0,\h) circle (\p) node(21) {};
\draw[fill] (0.309*\r,0.951*\r,\h) circle (\p) node(22) {};
\draw[fill] (-0.809*\r,0.588*\r,\h) circle (\p) node(23) {};
\draw[fill] (-0.809*\r,-0.588*\r,\h) circle (\p) node(24) {};
\draw[fill] (0.309*\r,-0.951*\r,\h) circle (\p) node(25) {};

\draw[fill] (\r,0,-\h) circle (\p) node(01) {};
\draw[fill] (0.309*\r,0.951*\r,-\h) circle (\p) node(02) {};
\draw[fill] (-0.809*\r,0.588*\r,-\h) circle (\p) node(03) {};
\draw[fill] (-0.809*\r,-0.588*\r,-\h) circle (\p) node(04) {};
\draw[fill] (0.309*\r,-0.951*\r,-\h) circle (\p) node(05) {};

\draw[fill] (2*\r-0.309*\r,0.951*\r,0) circle (\p) node(16) {};
\draw[fill] (2*\r+0.809*\r,0.588*\r,0) circle (\p) node(17) {};
\draw[fill] (2*\r+0.809*\r,-0.588*\r,0) circle (\p) node(18) {};
\draw[fill] (2*\r-0.309*\r,-0.951*\r,0) circle (\p) node(19) {};

\draw[fill] (2*\r-0.309*\r,0.951*\r,\h) circle (\p) node(26) {};
\draw[fill] (2*\r+0.809*\r,0.588*\r,\h) circle (\p) node(27) {};
\draw[fill] (2*\r+0.809*\r,-0.588*\r,\h) circle (\p) node(28) {};
\draw[fill] (2*\r-0.309*\r,-0.951*\r,\h) circle (\p) node(29) {};

\draw[fill] (2*\r-0.309*\r,0.951*\r,-\h) circle (\p) node(06) {};
\draw[fill] (2*\r+0.809*\r,0.588*\r,-\h) circle (\p) node(07) {};
\draw[fill] (2*\r+0.809*\r,-0.588*\r,-\h) circle (\p) node(08) {};
\draw[fill] (2*\r-0.309*\r,-0.951*\r,-\h) circle (\p) node(09) {};

\draw[black,fill opacity=\op,fill=\gris] (01.center) -- (02.center) -- (03.center) -- (04.center) -- (05.center) -- (01.center);
\draw[black,fill opacity=\op,fill=\gris] (11.center) -- (12.center) -- (13.center) -- (14.center) -- (15.center) -- (11.center);
\draw[black,fill opacity=\op,fill=\gris] (21.center) -- (22.center) -- (23.center) -- (24.center) -- (25.center) -- (21.center);

\draw[black,fill opacity=\op,fill=\gris] (01.center) -- (06.center) -- (07.center) -- (08.center) -- (09.center) -- (01.center);
\draw[black,fill opacity=\op,fill=\gris] (11.center) -- (16.center) -- (17.center) -- (18.center) -- (19.center) -- (11.center);
\draw[black,fill opacity=\op,fill=\gris] (21.center) -- (26.center) -- (27.center) -- (28.center) -- (29.center) -- (21.center);

\draw[black,fill opacity=\oop,fill=\griss] (01.center) -- (02.center) -- (12.center) -- (11.center) -- (01.center);
\draw[black,fill opacity=\oop,fill=\griss] (21.center) -- (22.center) -- (12.center) -- (11.center) -- (21.center);
\draw[black,fill opacity=\oop,fill=\griss] (03.center) -- (02.center) -- (12.center) -- (13.center) -- (03.center);
\draw[black,fill opacity=\oop,fill=\griss] (23.center) -- (22.center) -- (12.center) -- (13.center) -- (23.center);
\draw[black,fill opacity=\oop,fill=\griss] (03.center) -- (04.center) -- (14.center) -- (13.center) -- (03.center);
\draw[black,fill opacity=\oop,fill=\griss] (23.center) -- (24.center) -- (14.center) -- (13.center) -- (23.center);
\draw[black] (25.center) -- (15.center) -- (05.center);

\draw[black,fill opacity=\oop,fill=\griss] (01.center) -- (06.center) -- (16.center) -- (11.center) -- (01.center);
\draw[black,fill opacity=\oop,fill=\griss] (21.center) -- (26.center) -- (16.center) -- (11.center) -- (21.center);
\draw[black,fill opacity=\oop,fill=\griss] (07.center) -- (06.center) -- (16.center) -- (17.center) -- (07.center);
\draw[black,fill opacity=\oop,fill=\griss] (27.center) -- (26.center) -- (16.center) -- (17.center) -- (27.center);
\draw[black,fill opacity=\oop,fill=\griss] (07.center) -- (08.center) -- (18.center) -- (17.center) -- (07.center);
\draw[black,fill opacity=\oop,fill=\griss] (27.center) -- (28.center) -- (18.center) -- (17.center) -- (27.center);
\draw[black] (29.center) -- (19.center) -- (09.center);

\draw [->,ultra thick,black!30!red] (09) edge (11) (11) edge (25);
\draw [->,ultra thick,blue] (02) edge (11) (11) edge (29);

\node (a) at (1.5,-1,-1.5) {\color{red} \Large\bfseries $a$};
\node (b) at (1.5,1,1.5) {\color{blue} \Large\bfseries $b$};
\node (x) at (2,0.2,0) {\Large\bfseries $x_0$};
\end{tikzpicture}
\end{center}
\caption{A part of the complex $X_5$, with the axes of $a$ and $b$.}
\label{fig:X5}
\end{figure}

Assume first that $m$ is odd, then according to Brady and McCammond (see~\cite{brady_mccammond_artinthree}), there is an interesting presentation of $I_2(m)$ given by $I_2(m) = \<a,b \st w_m(a,b) = w_m(b,a)\> = \<t,u \st t^m=u^2\>$, where $t=ab$ and $u=w_m(a,b)$, so the central quotient $G$ of $I_2(m)$ is isomorphic to $ \<t,u \st t^m=u^2\> / \<t^m=u^2\> \simeq \Z/m\Z \star \Z/2\Z$. Consider the action of $G$ on the Bass-Serre $(m,2)$-biregular tree $T$, and consider the regular $m$-gonal complex $T_m$ obtained from $T$ by replacing the star of each vertex with valency $m$ by a regular $m$-gon with side length $1$, where $t$ acts on the base $m$-gon $P$ by a rotation of angle $\frac{4\pi}{m}$. Let us denote $p=\f{m-1}{2}$. Note that $a=t^{-p}u$ and $b=ut^{-p}$, and $t^p$ acts on the base $m$-gon by a rotation of angle $\frac{4p\pi}{m}=\frac{-2\pi}{m}$. This way, the axes of $a$ and $b$ acting on $T_m$ intersect the boundary of the $m$-gon $P$ in consecutive sides. Let $e \in T_m$ denote the intersection of the axes of $a$ and $b$, it is also the unique vertex fixed by $u=w_m(a,b)$ (see Figure~\ref{fig:X5}).

\begin{figure}
\begin{center}
\begin{tikzpicture}
\def \p {0.05}
\def \op {1}
\def \gris {black!10}
\draw[fill] (-1,0)+(0:1) circle (\p) node(1) {};
\draw[fill] (-1,0)+(72:1) circle (\p) node(2){} ;
\draw[fill] (-1,0)+(144:1) circle (\p) node(3) {};
\draw[fill] (-1,0)+(-144:1) circle (\p) node(4) {};
\draw[fill] (-1,0)+(-72:1) circle (\p) node(5) {};
\draw[fill] (1,0)+(180+72:1) circle (\p) node(6) {};
\draw[fill] (1,0)+(180+144:1) circle (\p) node(7) {};
\draw[fill] (1,0)+(180-144:1) circle (\p) node(8) {};
\draw[fill] (1,0)+(180-72:1) circle (\p) node(9) {};

\draw[draw=white] (125:1)+(2) circle (\p) node(10) {};
\draw[fill] (10)+(-55+72:1) circle (\p) node(11) {};
\draw[fill] (10)+(-55+144:1) circle (\p) node(12) {};
\draw[fill] (10)+(-55-144:1) circle (\p) node(13) {};
\draw[fill] (10)+(-55-72:1) circle (\p) node(14) {};
\draw[draw=white] (55:1)+(9) circle (\p) node(15) {};
\draw[fill] (15)+(-125+72:1) circle (\p) node(16) {};
\draw[fill] (15)+(-125+144:1) circle (\p) node(17) {};
\draw[fill] (15)+(-125-144:1) circle (\p) node(18) {};
\draw[fill] (15)+(-125-72:1) circle (\p) node(19) {};

\draw[draw=white] (-125:1)+(5) circle (\p) node(20) {};
\draw[fill] (20)+(55+72:1) circle (\p) node(21) {};
\draw[fill] (20)+(55+144:1) circle (\p) node(22) {};
\draw[fill] (20)+(55-144:1) circle (\p) node(23) {};
\draw[fill] (20)+(55-72:1) circle (\p) node(24) {};
\draw[draw=white] (-55:1)+(6) circle (\p) node(25) {};
\draw[fill] (25)+(125+72:1) circle (\p) node(26) {};
\draw[fill] (25)+(125+144:1) circle (\p) node(27) {};
\draw[fill] (25)+(125-144:1) circle (\p) node(28) {};
\draw[fill] (25)+(125-72:1) circle (\p) node(29) {};

\draw[black,fill opacity=\op,fill=\gris] (1.center) -- (2.center) -- (3.center) -- (4.center) -- (5.center) -- (1.center);
\draw[black,fill opacity=\op,fill=\gris] (1.center) -- (9.center) -- (8.center) -- (7.center) -- (6.center) -- (1.center);
\draw[black,fill opacity=\op,fill=\gris] (2.center) -- (11.center) -- (12.center) -- (13.center) -- (14.center) -- (2.center);
\draw[black,fill opacity=\op,fill=\gris] (9.center) -- (16.center) -- (17.center) -- (18.center) -- (19.center) -- (9.center);
\draw[black,fill opacity=\op,fill=\gris] (5.center) -- (21.center) -- (22.center) -- (23.center) -- (24.center) -- (5.center);
\draw[black,fill opacity=\op,fill=\gris] (6.center) -- (26.center) -- (27.center) -- (28.center) -- (29.center) -- (6.center);

\node (31) at ([shift=(120:0.5)]3) {};
\node (32) at ([shift=(120+108:0.5)]3) {};
\node (33) at ([shift=(-120:0.5)]4) {};
\node (34) at ([shift=(-120-108:0.5)]4) {};
\draw[black,fill opacity=\op,fill=\gris] (32.center) -- (3.center) -- (31.center);
\draw[black,fill opacity=\op,fill=\gris] (33.center) -- (4.center) -- (34.center);

\node (41) at ([shift=(60:0.5)]8) {};
\node (42) at ([shift=(60-108:0.5)]8) {};
\node (43) at ([shift=(-60:0.5)]7) {};
\node (44) at ([shift=(-60+108:0.5)]7) {};
\draw[black,fill opacity=\op,fill=\gris] (42.center) -- (8.center) -- (41.center);
\draw[black,fill opacity=\op,fill=\gris] (43.center) -- (7.center) -- (44.center);

\draw [->,ultra thick,black!30!red] (26) edge (6) (6) edge (1) (1) edge (5) (5) edge (24);
\draw [<-,ultra thick,blue] (19) edge (9) (9) edge (1) (1) edge (2) (2) edge (11);

\node (P) at ([xshift=-1cm]1) {\bfseries $P$};
\node (e) at ([xshift=-0.3cm]1) {\bfseries $e$};
\node (a) at ([yshift=-1cm]1) {\color{red} \Large\bfseries $a$};
\node (b) at ([yshift=+1cm]1) {\color{blue} \Large\bfseries $b$};

\end{tikzpicture}
\end{center}
\caption{A part of the complex $T_5$, with the axes of $a$ and $b$.}
\label{fig:T5}
\end{figure}
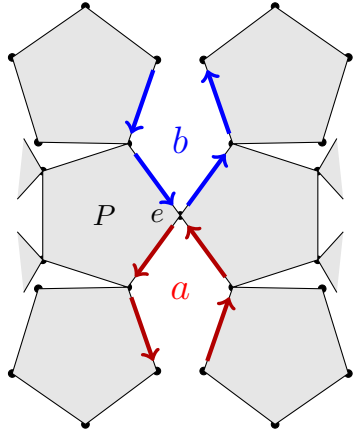

Consider the action of $I_2(m)$ on $\R$ by $a \cdot x = b \cdot x = x + \alpha$. We can endow $\R$ with the piecewise Euclidean simplicial structure where the vertex set is $\alpha \Z$.

Let $Y_m = T_m \times \R$, endowed with the diagonal action of $I_2(m)$, with basepoint $y_0=(e,0)$. The stabilizers of the points of $T_m$ are conjugated to either the cyclic subgroup spanned by $u$ or by $t$, and these subgroups act freely properly by translations on $\R$, we deduce that the action of $I_2(m)$ on $Y_m$ is free, with compact quotient. More precisely, since $I_2(m)$ acts transitively on vertices of $T_m$ and since the stabilizer of each vertex has $m$ orbits of vertices in $\R$, we deduce that $Y_m$ has exactly $m$ orbits of vertices. Furthermore, since $Y_m$ is locally finite, we deduce that $X_m = I_2(m) \bs Y_m$ is a compact locally CAT(0) space such that $\pi_1(X_m,x_0)$ is isomorphic to $I_2(m)$, where $x_0=I_2(m) \cdot y_0$ (see Figure~\ref{fig:X5}).

Let $Y_m^a$, $Y_m^b$ denote the axes of $a$ and $b$ through $y_0$, their image in $X_m$ define locally geodesic oriented loops $X_m^a$ and $X_m^b$, such that the angle at $x_0$ between $a^+$ and $b^+$(and similarly between $a^-$ and $b^-$) is equal to $\pi-2\arctan(\alpha)>\frac{4\pi}{5}$. 

Concerning the angle between $a^+$ and $b^-$, note that the $\R$ components of these vectors have opposite signs, hence the angle between $a^+$ and $b^-$ is strictly bigger than the angle between their $T_m$ components, which is precisely $\f{(m-2)\pi}{m}$. Hence
$$\sphericalangle_{x_0}(a^+,b^-) > \f{(m-2)\pi}{m} \geq \f{3\pi}{5},$$
and similarly between $a^-$ and $b^+$.

Up to rescaling $X_m$, since $X_m^a$ and $X_m^b$ have the same length, we can assume that they both have length $1$. 

\mk

Assume now that $m=2p$ is even, then according to Brady and McCammond (see~\cite{brady_mccammond_artinthree}), there is an interesting presentation of $I_2(m)$ given by $I_2(m) = \<a,b \st w_m(a,b) = w_m(b,a)\> = \<a,t \st at^p=t^pa\>$, where $t=ab$. In particular, $I_2(m)$ can be seen as the HNN extension of the group $\<t\> \simeq \Z$ with the subgroup $\<t^p\>$ and the identity map, with stable letter $a$.

Consider the action of $I_2(2p)$ on the Bass-Serre oriented $2p$-regular tree $T$. Let $T'$ denote the barycentric subdivision of $T$, it is an oriented $(2p,2)$-biregular tree. Consider the regular $2p$-gonal complex $T_{2p}$ obtained from $T'$ by replacing the star of each vertex with degree $2p$ by a regular $2p$-gon with side length $1$, such that $t$ acts on the base $2p$-gon $P$ by a rotation of angle $\frac{4\pi}{2p}$. The action of $\<t\>$ on the vertices of $P$ has two orbits, corresponding to the two possible orientations of edges adjacent to the base vertex.

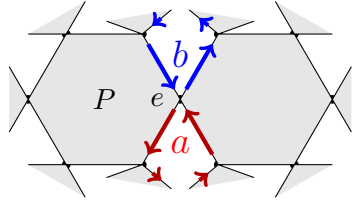
\begin{figure}
\begin{center}
\begin{tikzpicture}
\def \p {0.05}
\def \op {1}
\def \gris {black!10}
\draw[fill] (-1,0)+(0:1) circle (\p) node(1) {};
\draw[fill] (-1,0)+(60:1) circle (\p) node(2){} ;
\draw[fill] (-1,0)+(120:1) circle (\p) node(3) {};
\draw[fill] (-1,0)+(180:1) circle (\p) node(4) {};
\draw[fill] (-1,0)+(-120:1) circle (\p) node(5) {};
\draw[fill] (-1,0)+(-60:1) circle (\p) node(6) {};
\draw[fill] (1,0)+(180+60:1) circle (\p) node(7) {};
\draw[fill] (1,0)+(180+120:1) circle (\p) node(8) {};
\draw[fill] (1,0)+(180+180:1) circle (\p) node(9) {};
\draw[fill] (1,0)+(180-120:1) circle (\p) node(10) {};
\draw[fill] (1,0)+(180-60:1) circle (\p) node(11) {};

\draw[black,fill opacity=\op,fill=\gris] (1.center) -- (2.center) -- (3.center) -- (4.center) -- (5.center) -- (6.center) -- (1.center);
\draw[black,fill opacity=\op,fill=\gris] (1.center) -- (11.center) -- (10.center) -- (9.center) -- (8.center) -- (7.center) -- (1.center);

\node (2a) at ([shift=(40:0.5)]2) {};
\node (2b) at ([shift=(160:0.5)]2) {};
\draw[black,fill opacity=\op,fill=\gris] (2a.center) -- (2.center) -- (2b.center);
\node (3a) at ([shift=(60:0.5)]3) {};
\node (3b) at ([shift=(180:0.5)]3) {};
\draw[black,fill opacity=\op,fill=\gris] (3a.center) -- (3.center) -- (3b.center);
\node (4a) at ([shift=(120:0.5)]4) {};
\node (4b) at ([shift=(-120:0.5)]4) {};
\draw[black,fill opacity=\op,fill=\gris] (4a.center) -- (4.center) -- (4b.center);
\node (6a) at ([shift=(-40:0.5)]6) {};
\node (6b) at ([shift=(-160:0.5)]6) {};
\draw[black,fill opacity=\op,fill=\gris] (6a.center) -- (6.center) -- (6b.center);
\node (5a) at ([shift=(-60:0.5)]5) {};
\node (5b) at ([shift=(-180:0.5)]5) {};
\draw[black,fill opacity=\op,fill=\gris] (5a.center) -- (5.center) -- (5b.center);

\node (11a) at ([shift=(180-40:0.5)]11) {};
\node (11b) at ([shift=(180-160:0.5)]11) {};
\draw[black,fill opacity=\op,fill=\gris] (11a.center) -- (11.center) -- (11b.center);
\node (10a) at ([shift=(180-60:0.5)]10) {};
\node (10b) at ([shift=(180-180:0.5)]10) {};
\draw[black,fill opacity=\op,fill=\gris] (10a.center) -- (10.center) -- (10b.center);
\node (9a) at ([shift=(180-120:0.5)]9) {};
\node (9b) at ([shift=(180+120:0.5)]9) {};
\draw[black,fill opacity=\op,fill=\gris] (9a.center) -- (9.center) -- (9b.center);
\node (8a) at ([shift=(-120:0.5)]8) {};
\node (8b) at ([shift=(0:0.5)]8) {};
\draw[black,fill opacity=\op,fill=\gris] (8a.center) -- (8.center) -- (8b.center);
\node (7a) at ([shift=(-140:0.5)]7) {};
\node (7b) at ([shift=(-20:0.5)]7) {};
\draw[black,fill opacity=\op,fill=\gris] (7a.center) -- (7.center) -- (7b.center);

\draw [->,ultra thick,black!30!red] (7a) edge (7) (7) edge (1) (1) edge (6) (6) edge (6a);
\draw [->,ultra thick,blue] (2a) edge (2) (2) edge (1) (1) edge (11) (11) edge (11a);

\node (P) at ([xshift=-1cm]1) {\bfseries $P$};
\node (e) at ([xshift=-0.3cm]1) {\bfseries $e$};
\node (a) at ([yshift=-0.6cm]1) {\color{red} \Large\bfseries $a$};
\node (b) at ([yshift=+0.6cm]1) {\color{blue} \Large\bfseries $b$};

\end{tikzpicture}
\end{center}
\caption{A part of the tree $T_6$, with the axes of $a$ and $b$.}
\label{fig:T6}
\end{figure}

Since $b=a^{-1}t$, the axes of $a$ and $b$ acting on $T_{2p}$ intersect the boundary of the $2p$-gon $P$ in consecutive sides. Let $e \in T_m$ denote the intersection of the axes of $a$ and $b$ (see Figure~\ref{fig:T6}).

Consider the action of $I_2(m)$ on $\R$ by $a \cdot x = b \cdot x = x + \alpha$. We can endow $\R$ with the piecewise Euclidean simplicial structure where the vertex set is $\alpha \Z$.

Let $Y_m = T_m \times \R$, endowed with the diagonal action of $I_2(m)$, with basepoint $y_0=(e,0)$. The stabilizers of the vertices of $T_m$ are all equal to the kernel $Z(I_2(m))=\<t^p\>$ of the action. This cyclic subgroup acts freely properly by translations on $\R$, we deduce that the action of $I_2(m)$ on $Y_m$ is free. More precisely, since $I_2(m)$ acts transitively on vertices of $T_m$ and since the stabilizer of each vertex has $2p$ orbits of vertices in $\R$, we deduce that $Y_m$ has exactly $2p$ orbits of vertices. Furthermore, since $Y_m$ is locally finite, we deduce that $X_m = I_2(m) \bs Y_m$ is a compact locally CAT(0) space such that $\pi_1(X_m,x_0)$ is isomorphic to $I_2(m)$, where $x_0=I_2(m) \cdot y_0$.

Let $Y_m^a$, $Y_m^b$ denote the axes of $a$ and $b$ through $y_0$, their images in $X_m$ define locally geodesic oriented loops $X_m^a$ and $X_m^b$, such that the angles at $x_0$ between $a^+$, $b^+$ and between $a^-$, $b^-$ are both equal to
$$ \pi - 2\arctan(\alpha)>\frac{4\pi}{5}.$$

Concerning the angle between $a^+$ and $b^-$, note that the $\R$ components of these vectors have opposite signs, hence the angle between $a^+$ and $b^-$ is strictly bigger than the angle between their $T_m$ components, which is precisely $\f{(m-2)\pi}{m}$. Hence
$$\sphericalangle_{x_0}(a^+,b^-) > \f{(m-2)\pi}{m} \geq \f{3\pi}{5},$$
and similarly between $a^-$ and $b^+$. And, if $m \geq 6$, this angle is bigger than $\frac{2\pi}{3}$.

Up to rescaling $X_m$, since $X_m^a$ and $X_m^b$ have the same length, we can assume that they both have length $1$. We can also assume that, up to refining the piecewise Euclidean structure of $X_m$, the axes $X_m^a$ and $X_m^b$ lie in the $1$-skeleton.

\ep

\section{The general case of XXL type Artin groups}

We now describe a metric model for XXL type Artin groups, obtained by gluing the complexes obtained by dihedral Artin groups. This is the first part of Theorem~\ref{thm:main}.

\bthm
For every XXL type Artin group $A$, there exists a compact locally CAT(0) $3$-dimensional piecewise Euclidean complex $X_A$ and $x_0 \in X_A$ such that $\pi_1(X_A,x_0) \simeq A$.
\ethm

\bp
For each $s \in S$, let $X_s$ denote a circle with length $1$ and basepoint $x_0 \in X_s$, such that $\pi_1(X_s,x_0)$ will be identified with $\<s\>$. Let $E$ denote the set of all edges of $\Gamma$. For each $I \in E$, let $X_I$ denote a copy of $X_m$, where $m$ is the label of the edge $I$. Consider the following space
$$ X_A = \left(\bigcup_{I \in E} X_I \cup \bigcup_{s \in S} X_s \right)/ \sim,$$
where the identifications are given, for all $s \in S$ and $I=\{s,t\} \in {\cal S}_2$, by $X_s \sim X_{s,t}^s$. According to the Van Kampen Theorem, the fundamental group $\pi_1(X_A,x_0)$ is isomorphic to $A$. Up to refining the cell structure, we can assume that $X_A$ is a piecewise Euclidean cell complex $X_A$. In order to prove that $X_A$ is locally CAT(0), according to Gromov's link condition, it is sufficient to prove that the link of every vertex is CAT(1). For every edge $e$ of $X_A$, the link of $e$ in $X_A$ is the disjoint union of links of $e$ in all $X_I$'s that contain $e$. Since each $X_I$ is CAT(0), the link of $e$ in $X_A$ is CAT(1).

In other words, it is enough to prove that the link of every vertex of $X_A$ is large, i.e. every closed locally geodesic loop has length at least $2\pi$. Fix a vertex $x \in X_A$, and assume that $\ell$ is a locally geodesic loop in the link of $x$. We will prove that $\ell$ has length at least $2\pi$.

\mk

Assume first that $\ell$ is contained in a unique $X_{ab}$. Since $X_{ab}$ is CAT(0), the link of $x$ in $X_A$ is large, so $\ell$ has length at least $2\pi$.

\mk

Assume now that $\ell$ is contained in $X_{ab} \cup X_{bc}$. Since $X_b$ is convex in both $X_{ab}$ and $X_{bc}$, we know that $X_{ab} \cup X_{bc}$ is CAT(0), hence $\ell$ has length at least $2\pi$.

\mk

Assume now that $\ell$ is contained in $X_{ab} \cup X_{bc} \cup X_{ac}$, but not less than three. Recall that $\ell$ is locally geodesic and $X_a$ and $X_b$ are convex in $X_{ab}$. Therefore if $\ell$ enters $X_{ab}$ through $X_a$, it exits $X_{ab}$ through $X_b$, and similarly for $X_{bc}$ and $X_{ac}$. In particular, the length of $\ell$ is at least 
$$\sphericalangle_{x_0}(a^+,b^-) + \sphericalangle_{x_0}(b^-,c^+) + \sphericalangle_{x_0}(c^+,a^+) = 2 \times \frac{3\pi}{5} + \frac{4\pi}{5} =2\pi \mbox{ or }$$
$$\sphericalangle_{x_0}(a^+,b^+) + \sphericalangle_{x_0}(b^+,c^+) + \sphericalangle_{x_0}(c^+,a^+) = 3 \times \frac{4\pi}{5} > 2\pi.$$

\mk

Assume now that $\ell$ is contained in no fewer than four $X_I$'s. Then its length is at least $4 \times \frac{3\pi}{5} > 2\pi$.

\mk

In conclusion, every locally geodesic loop in the link of $x$ has length at least $2\pi$. So the link of $x$ is CAT(1), and $X_A$ is locally CAT(0).

\ep

Note that this construction is not sharp, meaning that we could also build this way a locally CAT(0) model for some Artin groups which are not of XXL type. However, the precise combinatorial conditions would not be very elegant to write down. Furthermore, such a construction cannot be adapted to take into account the $(3,3,3)$ triangle Artin group for instance, which is known by Brady and McCammond (see~\cite{brady_mccammond_artinthree}) to be CAT(0) using another complex.

\section{A rank one geodesic}

We will now prove that the locally CAT(0) complex we built for XXL type Artin groups has rank $1$, meaning that there exists a periodic geodesic in the universal cover which does not bound any flat half-plane. Fix $\alpha \in (0,\tan (\f{\pi}{10}))$. We start by looking at a specific loop for the complex for dihedral Artin groups, first in the odd case.

\blem \label{lem:ell_odd} For every odd $m \geq 5$, there exists a locally geodesic oriented simple loop $\ell$ in $X_m$ based at $x_0$ such that, if we denote $\ell^+,\ell^- \in \lk_{x_0}(X_m)$ the images in the link of $x_0$ of the positive and negative sides of the loop $\ell$, we have:
\bit
\item $X_m^a \cap \ell = X_m^b \cap \ell = \{x_0\}$,
\item $\sphericalangle_{x_0}(a^+,\ell^+),\sphericalangle_{x_0}(a^-,\ell^-) > \frac{2\pi}{5}$.
\item $\sphericalangle_{x_0}(b^-,\ell^+),\sphericalangle_{x_0}(b^+,\ell^-) > \frac{\pi}{5}$.
\item $\sphericalangle_{x_0}(a^-,\ell^+),\sphericalangle_{x_0}(a^+,\ell^-),\sphericalangle_{x_0}(b^+,\ell^+),\sphericalangle_{x_0}(b^-,\ell^-) > \frac{4\pi}{5}$.
\eit
\elem

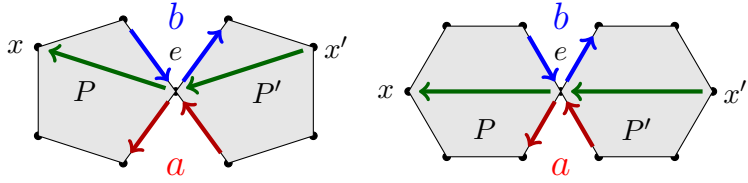
\begin{figure}
\begin{center}
\begin{tikzpicture}
\def \p {0.05}
\def \op {1}
\draw[fill] (-1,0)+(0:1) circle (\p) node(1) {};
\draw[fill] (-1,0)+(72:1) circle (\p) node(2){} ;
\draw[fill] (-1,0)+(144:1) circle (\p) node(3) {};
\draw[fill] (-1,0)+(-144:1) circle (\p) node(4) {};
\draw[fill] (-1,0)+(-72:1) circle (\p) node(5) {};
\draw[fill] (1,0)+(180+72:1) circle (\p) node(6) {};
\draw[fill] (1,0)+(180+144:1) circle (\p) node(7) {};
\draw[fill] (1,0)+(180-144:1) circle (\p) node(8) {};
\draw[fill] (1,0)+(180-72:1) circle (\p) node(9) {};

\draw[black,fill opacity=\op,fill=black!10] (1.center) -- (2.center) -- (3.center) -- (4.center) -- (5.center) -- (1.center);
\draw[black,fill opacity=\op,fill=black!10] (1.center) -- (9.center) -- (8.center) -- (7.center) -- (6.center) -- (1.center);

\draw [->,ultra thick,black!30!red] (6) edge (1) (1) edge (5);
\draw [<-,ultra thick,blue] (9) edge (1) (1) edge (2);
\draw [<-,ultra thick,black!60!green] (3) edge (1) (1) edge (8);

\node (P) at ([xshift=-1.2cm]1) {\bfseries $P$};
\node (P') at ([xshift=1.2cm]1) {\bfseries $P'$};
\node (e) at ([yshift=0.5cm]1) {\bfseries $e$};
\node (x) at ([xshift=-0.3cm]3) {\bfseries $x$};
\node (x') at ([xshift=0.3cm]8) {\bfseries $x'$};
\node (a) at ([yshift=-1cm]1) {\color{red} \Large\bfseries $a$};
\node (b) at ([yshift=+1cm]1) {\color{blue} \Large\bfseries $b$};

\end{tikzpicture}
\begin{tikzpicture}
\def \p {0.05}
\def \op {1}
\draw[fill] (-1,0)+(0:1) circle (\p) node(1) {};
\draw[fill] (-1,0)+(60:1) circle (\p) node(2){} ;
\draw[fill] (-1,0)+(120:1) circle (\p) node(3) {};
\draw[fill] (-1,0)+(180:1) circle (\p) node(4) {};
\draw[fill] (-1,0)+(-120:1) circle (\p) node(5) {};
\draw[fill] (-1,0)+(-60:1) circle (\p) node(6) {};
\draw[fill] (1,0)+(180+60:1) circle (\p) node(7) {};
\draw[fill] (1,0)+(180+120:1) circle (\p) node(8) {};
\draw[fill] (1,0)+(180+180:1) circle (\p) node(9) {};
\draw[fill] (1,0)+(180-120:1) circle (\p) node(10) {};
\draw[fill] (1,0)+(180-60:1) circle (\p) node(11) {};

\draw[black,fill opacity=\op,fill=black!10] (1.center) -- (2.center) -- (3.center) -- (4.center) -- (5.center) -- (6.center) -- (1.center);
\draw[black,fill opacity=\op,fill=black!10] (1.center) -- (11.center) -- (10.center) -- (9.center) -- (8.center) -- (7.center) -- (1.center);

\draw [->,ultra thick,black!30!red] (7) edge (1) (1) edge (6);
\draw [<-,ultra thick,blue] (11) edge (1) (1) edge (2);
\draw [<-,ultra thick,black!60!green] (4) edge (1) (1) edge (9);

\node (P) at ([xshift=-1cm,yshift=-0.5cm]1) {\bfseries $P$};
\node (P') at ([xshift=1cm,yshift=-0.5cm]1) {\bfseries $P'$};
\node (e) at ([yshift=0.5cm]1) {\bfseries $e$};
\node (x) at ([xshift=-0.3cm]4) {\bfseries $x$};
\node (x') at ([xshift=0.3cm]9) {\bfseries $x'$};
\node (a) at ([yshift=-1cm]1) {\color{red} \Large\bfseries $a$};
\node (b) at ([yshift=+1cm]1) {\color{blue} \Large\bfseries $b$};

\end{tikzpicture}
\end{center}
\caption{The construction of the loop in $X_5$ and in $X_6$}
\label{fig:ell}
\end{figure}

\bp
In the polygonal complex $T_m$, consider the two $m$-gons $P,P'$ adjacent to the base vertex $e$. Consider the vertex $x \in P$ such that $e$ and $x$ are "almost opposite" in $P$, i.e. form an angle of $\f{2p\pi}{m}$ from the center of $P$, where $m=2p+1$. Consider the unique vertex $x' \in P'$ such that $(x,0)$ and $(x',0)$ are in the same $I_2(m)$-orbit in $Y_m$, then $x'$ and $e$ form an angle of $\f{-2p\pi}{m}$ from the center of $P'$ (see Figure~\ref{fig:ell}).

The piecewise Euclidean path from $x'$ to $x$ consisting of the two segments from $x'$ to $e$ and from $e$ to $x$ projects to a locally geodesic oriented simple loop $\ell$ in the image of $T_m \times \{0\}$ in $X_m$. By construction, we have $X_m^a \cap \ell = X_m^b \cap \ell = \{x_0\}$. Furthermore, we have
\bit
\item $\sphericalangle_{x_0}(a^+,\ell^+)=\sphericalangle_{x_0}(a^-,\ell^-) > \f{p\pi}{m} \geq \frac{2\pi}{5}$ and
\item $\sphericalangle_{x_0}(b^-,\ell^+)=\sphericalangle_{x_0}(b^+,\ell^-) > \f{(p-1)\pi}{m} \geq \frac{\pi}{5}$.
\eit

Furthermore, since the angle at $e$ inside $T_m$ between $a^-$ and $\ell^+$ is infinite, we deduce that in the whole space $Y_m=T_m \times \R$ we have $\sphericalangle_{x_0}(a^-,\ell^+)=\pi-\arctan(\alpha)>\f{9\pi}{10}>\f{4\pi}{5}$. Similarly
$$\sphericalangle_{x_0}(a^-,\ell^+)=\sphericalangle_{x_0}(a^+,\ell^-)=\sphericalangle_{x_0}(b^+,\ell^+)=\sphericalangle_{x_0}(b^-,\ell^-) > \frac{4\pi}{5}.$$
\ep

We now turn to the even dihedral Artin groups.

\blem \label{lem:ell_even} For every even $m \geq 6$, there exists a locally geodesic oriented simple loop $\ell$ in $X_m$ based at $x_0$ such that, if we denote $\ell^+,\ell^- \in \lk_{x_0}(X_m)$ the images in the link of $x_0$ of the positive and negative sides of the loop $\ell$, we have:
\bit
\item $X_m^a \cap \ell = X_m^b \cap \ell = \{x_0\}$,
\item $\sphericalangle_{x_0}(a^+,\ell^+),\sphericalangle_{x_0}(a^-,\ell^-) > \frac{\pi}{3}$.
\item $\sphericalangle_{x_0}(b^-,\ell^+),\sphericalangle_{x_0}(b^+,\ell^-) > \frac{\pi}{3}$.
\item $\sphericalangle_{x_0}(a^-,\ell^+),\sphericalangle_{x_0}(a^+,\ell^-),\sphericalangle_{x_0}(b^+,\ell^+),\sphericalangle_{x_0}(b^-,\ell^-) > \frac{4\pi}{5}$.
\eit
\elem

\bp
In the polygonal complex $T_m$, consider the two $m$-gons $P,P'$ adjacent to the base vertex $e$. Consider the vertex $x \in P$ such that $e$ and $x$ are opposite in $P$. Consider the unique vertex $x' \in P'$ such that $(x,0)$ and $(x',0)$ are in the same $I_2(m)$-orbit in $Y_m$, then $x'$ and $e$ are opposite in $P'$ (see Figure~\ref{fig:ell}).

The piecewise Euclidean path from $x'$ to $x$ consisting of the two segments from $x'$ to $e$ and from $e$ to $x$ projects to a locally geodesic oriented simple loop $\ell$ in the image of $T_m \times \{0\}$ in $X_m$. By construction, we have $X_m^a \cap \ell = X_m^b \cap \ell = \{x_0\}$. Furthermore, we have
$$\sphericalangle_{x_0}(a^+,\ell^+)=\sphericalangle_{x_0}(a^-,\ell^-) = \sphericalangle_{x_0}(b^-,\ell^+)=\sphericalangle_{x_0}(b^+,\ell^-) > \f{(m-2)\pi}{2m} \geq \frac{\pi}{3}.$$

Furthermore, since the angle at $e$ inside $T_m$ between $a^-$ and $\ell^+$ is infinite, we deduce that in the whole space $Y_m=T_m \times \R$ we have $\sphericalangle_{x_0}(a^-,\ell^+)=\pi-\arctan(\alpha)>\f{9\pi}{10}>\f{4\pi}{5}$. Similarly
$$\sphericalangle_{x_0}(a^-,\ell^+)=\sphericalangle_{x_0}(a^+,\ell^-)=\sphericalangle_{x_0}(b^+,\ell^+)=\sphericalangle_{x_0}(b^-,\ell^-) > \frac{4\pi}{5}.$$
\ep

We can now prove that the complex $X_A$ has an isometry of rank one, thus proving the second part of Theorem~\ref{thm:main}.

\bthm Assume that $A(\Gamma)$ is an XXL type Artin group with at least three generators. Then there exists a locally geodesic loop in $X_A$ whose lifts in $\widetilde{X_A}$ have rank $1$, i.e. do not bound flat half-planes. \ethm

\bp
\bit
\item If $\Gamma$ has no edge, then $A(\Gamma)$ is the free group on $S$ and $X_A$ is a wedge of $|S|$ circles. So every geodesic in the tree $\widetilde{X_A}$ has rank $1$.

\item Assume now that $\Gamma$ has at least one edge labeled by some odd number $m \geq 5$, between $a$ and $b$. Fix $c \in S \bs \{a,b\}$. Let $\ell_{ab} \subset X_{ab} \subset X_A$ denote the loop given by Lemma~\ref{lem:ell_odd}, and consider the oriented loop $X_c \subset X_A$. Then consider the concatenation $\ell=\ell_{ab} \cdot X_c$. We will prove that the angle at $x_0$ between the incoming loop $c^-$ and the outgoing loop $\ell_{ab}^+$ is bigger than $\pi$.

\mk

By construction, in the link of $x_0$, every path from $\ell_{ab}^+$ to $c^-$ must pass through one of $\{a^+,a^-,b^+,b^-\}$. Let us compute the four quantities:
\beq
\sphericalangle_{x_0}(\ell_{ab}^+,a^+)+\sphericalangle_{x_0}(a^+,c^-) &>& \f{2\pi}{5} + \f{3\pi}{5}=\pi, \\
\sphericalangle_{x_0}(\ell_{ab}^+,a^-)+\sphericalangle_{x_0}(a^-,c^-) &>& \f{4\pi}{5} + \f{4\pi}{5}>\pi, \\
\sphericalangle_{x_0}(\ell_{ab}^+,b^+)+\sphericalangle_{x_0}(b^+,c^-) &>& \f{4\pi}{5} + \f{3\pi}{5}>\pi, \\
\sphericalangle_{x_0}(\ell_{ab}^+,b^-)+\sphericalangle_{x_0}(b^-,c^-) &>& \f{\pi}{5} + \f{4\pi}{5}=\pi.
\eeq

We deduce that the distance in the link of $x_0$ between $\ell_{ab}^+$ and $c^-$ is bigger than $\pi$. Similarly, the distance in the link of $x_0$ between $\ell_{ab}^-$ and $c^+$ is also bigger than $\pi$. 

\item Assume now that $\Gamma$ has at least one edge, and that all edges are labeled by even numbers. Consider some edge between $a$ and $b$, labeled by some even $m \geq 6$. Fix $c \in S \bs \{a,b\}$. Let $\ell_{ab} \subset X_{ab} \subset X_A$ denote the loop given by Lemma~\ref{lem:ell_even}, and consider the oriented loop $X_c \subset X_A$. Then consider the concatenation $\ell=\ell_{ab} \cdot X_c$. We will prove that the angle at $x_0$ between the incoming loop $c^-$ and the outgoing loop $\ell_{ab}^+$ is bigger than $\pi$.

\mk

By construction, in the link of $x_0$, every path from $\ell_{ab}^+$ to $c^-$ must pass through one of $\{a^+,a^-,b^+,b^-\}$. Let us compute the four quantities:
\beq
\sphericalangle_{x_0}(\ell_{ab}^+,a^+)+\sphericalangle_{x_0}(a^+,c^-) &>& \f{\pi}{3} + \f{2\pi}{3}=\pi, \\
\sphericalangle_{x_0}(\ell_{ab}^+,a^-)+\sphericalangle_{x_0}(a^-,c^-) &>& \f{4\pi}{5} + \f{4\pi}{5}>\pi, \\
\sphericalangle_{x_0}(\ell_{ab}^+,b^+)+\sphericalangle_{x_0}(b^+,c^-) &>& \f{4\pi}{5} + \f{2\pi}{3}>\pi, \\
\sphericalangle_{x_0}(\ell_{ab}^+,b^-)+\sphericalangle_{x_0}(b^-,c^-) &>& \f{\pi}{3} + \f{4\pi}{5}>\pi.
\eeq

We deduce that the distance in the link of $x_0$ between $\ell_{ab}^+$ and $c^-$ is bigger than $\pi$. Similarly, the distance in the link of $x_0$ between $\ell_{ab}^-$ and $c^+$ is also bigger than $\pi$. 
\eit

\mk

In conclusion, in each of the last two cases, $\ell$ is a locally geodesic loop in $X_A$ such that the angle at (each of the two passings at) $x_0$ is bigger than $\pi$. In particular, any lift of $\ell$ in $\widetilde{X_A}$ does not bound a flat half-plane, so it has rank $1$.
\ep

\sign

\bibliographystyle{../../smfalpha_perso}
\bibliography{../../bibli}

\end{document}